# Dirichlet series and series with Stirling numbers

Khristo N. Boyadzhiev

Department of Mathematics
Ohio Northern University, Ada, Ohio 45810, USA
k-boyadzhiev@onu.edu

**Abstract**

This paper presents a number of identities for Dirichlet series and series with Stirling numbers of the first kind. As coefficients for the Dirichlet series we use Cauchy numbers of the first and second kinds, hyperharmonic numbers, derangement numbers, binomial coefficients, central binomial coefficients, and Catalan numbers.

## 1. Introduction

We consider the unsigned Stirling numbers of the first kind $\begin{bmatrix} n \\ k \end{bmatrix}$ defined by the equation

$$(1) \qquad \frac{1}{x^{k+1}} = \sum_{n=k}^{\infty} \begin{bmatrix} n \\ k \end{bmatrix} \frac{1}{x(x+1)...(x+n)}$$

(see [2], [15, p. 29 and p. 171]) or by the well-known exponential generating function [8]

$$(2) \qquad \frac{(-1)^k}{k!} \ln^k(1-x) = \sum_{n=k}^{\infty} \begin{bmatrix} n \\ k \end{bmatrix} \frac{x^n}{n!} \quad (|x|<1).$$

These numbers have a very strong presence in combinatorics and also in classical analysis. For example, it is known that for $k \geq 1$

$$(3) \qquad \zeta(k+1) = \sum_{n=k}^{\infty} \begin{bmatrix} n \\ k \end{bmatrix} \frac{1}{n!n}$$

where $\zeta(s)$ is Riemann's zeta function. This representation appears in Jordan's book [12] (see equation (6) on p. 166 and equation (11) on p. 194 in [12]; also the more general formula on p. 343 in this book). Interesting comments and a new proof can be found in Adamchik's paper [1].

The corresponding result for the Hurwitz zeta function $\zeta(s,a)$ was proved in [2]

(4) $$\zeta(k+1,a) = \Gamma(a)\sum_{n=k}^{\infty}\begin{bmatrix}n\\k\end{bmatrix}\frac{1}{n\Gamma(n+a)} \quad (a>0, k\geq 1)$$

together with the representation

(5) $$\sum_{p=1}^{\infty}\frac{H_p}{p^{k+1}} = \sum_{n=k}^{\infty}\begin{bmatrix}n\\k\end{bmatrix}\frac{\psi'(n)}{n!} \quad (k\geq 1)$$

where $H_n = 1 + \frac{1}{2} + \ldots + \frac{1}{n}$, $H_0 = 0$, are the harmonic numbers and $\psi(x) = \frac{d}{dx}\ln\Gamma(x)$ is the digamma function. A representation for the polylogarithm $\mathrm{Li}_s(x) = \sum_{n=1}^{\infty}\frac{x^n}{n^s}$ in the spirit of (3) was also proved in [2]. Other results in this line were obtained by Wang and Chen [16].

In this paper we want to develop a method of replacing $\zeta(k+1)$ by other Dirichlet series in the spirit of (4) and (5) and obtaining representations in terms of Stirling numbers. More precisely, we consider the following problem arising from the above representations: given a Dirichlet series of the form

(6) $$F(s) = \sum_{n=1}^{\infty}\frac{a_{n-1}}{n^s}$$

we want to find the coefficients $b_n$ so that

$$F(k+1) = \sum_{n=k}^{\infty}\begin{bmatrix}n\\k\end{bmatrix}b_n.$$

In the following section (Section 2) we state our main theorem which gives a solution to this problem for the case when the numbers $a_k$ are the coefficients of an analytic function. Then in Section 3 we derive from this theorem a number of corollaries. For illustration, here are some identities proved in Section 3

$$k+1-\sum_{j=1}^{k}\zeta(j+1) = \sum_{n=k}^{\infty}\begin{bmatrix}n\\k\end{bmatrix}\frac{1}{n!(n+1)^2} \quad (k\geq 1)$$

$$\sum_{p=0}^{\infty}\frac{1}{p!(p+1)^{k+1}} = \sum_{n=k}^{\infty}\begin{bmatrix}n\\k\end{bmatrix}\left(e - \sum_{j=0}^{n}\frac{1}{j!}\right) \quad (k\geq 0)$$

$$\sum_{p=0}^{\infty}\frac{C_p}{4^p(p+1)^{k+1}} = 4\sum_{n=k}^{\infty}\begin{bmatrix}n\\k\end{bmatrix}\frac{\beta(2n+2)}{n!} \quad (k\geq 0)$$

$$\sum_{p=1}^{\infty}\frac{C_{p-1}}{4^p(p+1)^{k+1}}+\sum_{n=k}^{\infty}\begin{bmatrix}n\\k\end{bmatrix}\frac{1}{n!(2n+3)}=\frac{1}{2} \quad (k\geq 0)$$

where $C_p$ are the Catalan numbers and $\beta(x)$ is Nielsen's beta function.

## 2. Main results

**Theorem 1.** *Suppose the function*

$$f(x)=\sum_{k=0}^{\infty}a_k x^k$$

*is analytic on the unit disc $|x|<1$. Then we have the representations*

(7) $$\sum_{p=1}^{\infty}\frac{a_{p-1}}{p^{k+1}}=\sum_{p=0}^{\infty}\frac{a_p}{(p+1)^{k+1}}=\sum_{n=k}^{\infty}\begin{bmatrix}n\\k\end{bmatrix}\frac{1}{n!}\int_0^1 f(x)(1-x)^n\,dx$$

(8) $$\sum_{p=0}^{\infty}\frac{a_p}{(p+1)^{k+1}}=\frac{(-1)^k}{k!}\int_0^1 f(x)(\ln x)^k\,dx.$$

*Proof.* In the representation (1) let $x=p\geq 1$ be a positive integer and write

(9) $$\frac{1}{p^{k+1}}=\sum_{n=k}^{\infty}\begin{bmatrix}n\\k\end{bmatrix}\frac{1}{p(p+1)\ldots(p+n)}.$$

We know that

$$\frac{1}{p(p+1)\ldots(p+n)}=\frac{1}{n!}B(p,n+1)=\frac{1}{n!}\int_0^1 x^{p-1}(1-x)^n\,dx$$

where $B(p,q)$ is Euler's beta function defined by

$$B(p,q)=\frac{\Gamma(p)\Gamma(q)}{\Gamma(p+q)}=\int_0^1 x^{p-1}(1-x)^{q-1}\,dx.$$

Next, multiplying both sides of (9) by $a_{p-1}$, summing for $p=1,2,\ldots$, and exchanging the order of summation we have

$$\sum_{p=1}^{\infty}\frac{a_{p-1}}{p^{k+1}}=\sum_{n=k}^{\infty}\begin{bmatrix}n\\k\end{bmatrix}\sum_{p=1}^{\infty}\frac{a_{p-1}}{p(p+1)\ldots(p+n)}=\sum_{n=k}^{\infty}\begin{bmatrix}n\\k\end{bmatrix}\frac{1}{n!}\int_0^1 f(x)(1-x)^n\,dx.$$

This proves (7). The representation (8) follows from (7) by bringing the summation inside the integral and using equation (2). It can be proved also directly by expanding $f(x)$ and integrating term by term. The proof is completed.

**Some remarks**. A Dirichlet series of the form (6) is called convergent, if it is absolutely convergent for some $s_0$. In this case it is absolutely and uniformly convergent for $\operatorname{Re}(s) > s_0$ [9, 11]. Also, as shown by Charles Jordan in [12, p.161] the Stirling numbers of the first kind have asymptotic growth

$$(10) \qquad \begin{bmatrix} n \\ k \end{bmatrix} \sim \frac{(n-1)!}{(k-1)!}(\gamma + \ln n)^{k+1}$$

for $k \geq 1$ fixed. Here $\gamma$ is Euler's constant. We will keep this asymptotic in mind for the corollaries that follow.

Taking the function

$$f(x) = \frac{1}{1-x} = \sum_{n=0}^{\infty} x^n, \quad |x| < 1; \quad a_n = 1 \ (n = 0, 1, ...)$$

we have for $n \geq 1$

$$\int_0^1 f(x)(1-x)^n \, dx = \int_0^1 (1-x)^{n-1} dx = \frac{1}{n}$$

and (3) follows immediately from (7). Our main applications are given in the next section.

## 3. Corollaries

We first give a new proof to equation (5). Applying (7) to the function

$$f(x) = -\frac{\ln(1-x)}{x(1-x)} = \sum_{n=0}^{\infty} H_{n+1} x^n \quad (|x| < 1)$$

we find

$$\sum_{p=1}^{\infty} \frac{H_p}{p^{k+1}} = \sum_{n=k}^{\infty} \begin{bmatrix} n \\ k \end{bmatrix} \frac{1}{n!} \left\{ -\int_0^1 \frac{\ln(1-x)(1-x)^{n-1}}{x} dx \right\}.$$

With the substitution $1 - x = e^{-t}$ this integral becomes

$$\int_0^{\infty} \frac{te^{-nt}}{1-e^{-t}} dt = \sum_{m=0}^{\infty} \frac{1}{(m+n)^2} = \psi'(n)$$

and (5) follows. It is good to mention that

$$\psi'(n) = \frac{\pi^2}{6} - \left(1 + \frac{1}{2^2} + \ldots + \frac{1}{(n-1)^2}\right).$$

### 3.1. Dirichlet series with hyperharmonic numbers

Let $h_n^{(r)}$ be the hyperharmonic numbers defined by the equation

$$h_n^{(r+1)} = \binom{n+r}{r}(H_{n+r} - H_r)$$

for integers $n, r \geq 0$ [8]. When $r = 0$, $h_n^{(1)} = H_n$ are the harmonic numbers. The generating function for $h_n^{(r)}$ is given by

$$\sum_{n=0}^{\infty} h_n^{(r)} x^n = -\frac{\ln(1-x)}{(1-x)^r} \quad (|x| < 1).$$

Let $n > r$. We apply our theorem to the function $f(x) = -\frac{\ln(1-x)}{(1-x)^{r+1}}$ where the parameter $r+1$ is used for technical convenience. It is easy to compute

$$\int_0^1 f(x)(1-x)^n \, dx = -\int_0^1 \ln(1-x)(1-x)^{n-r-1} \, dx = \frac{1}{n-r} \int_0^1 \ln(1-x) \, d(1-x)^{n-r}$$

$$= \frac{1}{n-r} \ln(1-x)(1-x)^{n-r} \Big|_0^1 + \frac{1}{n-r} \int_0^1 (1-x)^{n-r-1} \, dx = \frac{1}{(n-r)^2}.$$

This gives the following.

**Corollary 2**. *For every $k > r \geq 0$ we have the identity*

(11) $$\sum_{p=1}^{\infty} \frac{h_{p-1}^{(r+1)}}{p^{k+1}} = \sum_{n=k}^{\infty} \begin{bmatrix} n \\ k \end{bmatrix} \frac{1}{n!(n-r)^2}.$$

In particular, for $r = 0$ we have (compare to (5))

(12) $$\sum_{p=1}^{\infty} \frac{H_{p-1}}{p^{k+1}} = \sum_{n=k}^{\infty} \begin{bmatrix} n \\ k \end{bmatrix} \frac{1}{n!n^2}.$$

### 3.2. On a result of Victor Adamchik

Adamchik in [1] discussed series of the form

$$G(k,q) = \sum_{n=k}^{\infty} \begin{bmatrix} n \\ k \end{bmatrix} \frac{1}{n!n^q}$$

and showed that

(13) $$G(k,q) = G(q,k)$$

that is,

$$\sum_{n=k}^{\infty} \begin{bmatrix} n \\ k \end{bmatrix} \frac{1}{n!n^q} = \sum_{n=q}^{\infty} \begin{bmatrix} n \\ q \end{bmatrix} \frac{1}{n!n^k}.$$

With $q=2$ this property implies equation (12) since $\begin{bmatrix} n \\ 2 \end{bmatrix} = (n-1)!H_{n-1}$.

In the next result we will connect two different series of Stirling numbers of the first kind.

**Corollary 3.** *For any two integers $q \geq 0, k \geq 0$ we have*

(14) $$\sum_{n=q}^{\infty} \begin{bmatrix} n \\ q \end{bmatrix} \frac{1}{n!(n+1)^{k+1}} = \sum_{n=k}^{\infty} \begin{bmatrix} n \\ k \end{bmatrix} \frac{1}{n!(n+1)^{q+1}}$$

*Proof.* We take the function

$$f(x) = \frac{(-1)^q}{q!} \ln^q(1-x) = \sum_{n=0}^{\infty} \begin{bmatrix} n \\ q \end{bmatrix} \frac{x^n}{n!}$$

where $q \geq 0$ is an integer. With

$$a_{p-1} = \frac{1}{(p-1)!} \begin{bmatrix} p-1 \\ q \end{bmatrix}$$

we find from the general formula (7)

$$\sum_{p=1}^{\infty} \begin{bmatrix} p-1 \\ q \end{bmatrix} \frac{1}{p^{k+1}(p-1)!} = \sum_{n=k}^{\infty} \begin{bmatrix} n \\ k \end{bmatrix} \frac{(-1)^q}{n!q!} \int_0^1 \ln^q(1-x)(1-x)^n dx$$

and with the substitution $1-x = e^{-t}$ we compute

$$\int_0^1 \ln^q(1-x)(1-x)^n dx = (-1)^q \int_0^{\infty} t^q e^{-(n+1)t} dt = \frac{(-1)^q q!}{(n+1)^{q+1}}$$

which gives the representation

$$\sum_{p=1}^{\infty} \begin{bmatrix} p-1 \\ q \end{bmatrix} \frac{1}{p^{k+1}(p-1)!} = \sum_{n=k}^{\infty} \begin{bmatrix} n \\ k \end{bmatrix} \frac{1}{n!(n+1)^{q+1}}.$$

Replacing $p$ by $n+1$ gives the desired result for all $q \geq 0, \ k \geq 0$, qed.

This shows a symmetry like the one in (13) for $G(k,q)$.

For $q=0$ we find for every $k \geq 0$

(15) $$\sum_{n=k}^{\infty} \begin{bmatrix} n \\ k \end{bmatrix} \frac{1}{(n+1)!} = 1.$$

Note that the same result follows from (1) for $x=1$. The value of this remarkable series is independent of $k$. For $k=0$ on the left hand side we have only one term, 1. For $k \geq 1$ fixed from (10) we find the asymptotic behavior for large $n$

$$\begin{bmatrix} n \\ k \end{bmatrix} \frac{1}{(n+1)!} \sim \frac{(\gamma + \ln n)^{k-1}}{n(n+1)(k-1)!}$$

which shows a very slow convergence. For instance,

$$\sum_{n=5}^{1000} \begin{bmatrix} n \\ 5 \end{bmatrix} \frac{1}{(n+1)!} < 0.8$$

For $q=1, k \geq 1$, we have $\begin{bmatrix} p \\ 1 \end{bmatrix} = (p-1)!$ and (7) implies the closed form evaluation

(16) $$\sum_{n=k}^{\infty} \begin{bmatrix} n \\ k \end{bmatrix} \frac{1}{n!(n+1)^2} = \sum_{p=1}^{\infty} \frac{1}{p(p+1)^{k+1}} = k+1 - \sum_{j=1}^{k} \zeta(j+1)$$

where the second equality follows from the recurrence relation

$$\sum_{p=1}^{\infty} \frac{1}{p(p+1)^{k+1}} = \sum_{p=1}^{\infty} \frac{1+p-p}{p(p+1)^{k+1}} = \sum_{p=1}^{\infty} \frac{1}{p(p+1)^k} - \zeta(k+1) + 1$$

$$= \sum_{p=1}^{\infty} \frac{1}{p(p+1)^{k-1}} - \zeta(k+1) - \zeta(k) + 2 \quad \text{etc.}$$

For $q=2$ we have $\begin{bmatrix} p \\ 2 \end{bmatrix} = (p-1)! H_{p-1}$ and therefore,

(17) $$\sum_{n=k}^{\infty} \begin{bmatrix} n \\ k \end{bmatrix} \frac{1}{n!(n+1)^3} = \sum_{p=1}^{\infty} \frac{H_{p-1}}{p(p+1)^{k+1}}$$

etc.

### 3.3. Dirichlet series with Cauchy numbers

The Cauchy numbers of the first kind $c_n$ and second kind $d_n$ are interesting combinatorial numbers. They are defined by the generating functions

$$\frac{x}{\ln(x+1)} = \sum_{n=0}^{\infty} \frac{c_n}{n!} x^n$$

$$\frac{-x}{(1-x)\ln(1-x)} = \sum_{n=0}^{\infty} \frac{d_n}{n!} x^n$$

where $|x|<1$ (see [7, p. 294] and [3, 14]). Now consider the function

$$f(x) = \frac{-x}{\ln(1-x)} = \sum_{n=0}^{\infty} \frac{(-1)^n c_n}{n!} x^n \quad (|x|<1).$$

Using again the substitution $1-x = e^{-t}$ we compute

$$\int_0^1 f(x)(1-x)^n dx = \int_0^{\infty} \frac{1-e^{-t}}{t} e^{-(n+1)t} dt = \int_0^{\infty} \frac{e^{-(n+1)t} - e^{-(n+2)t}}{t} dt = \ln\frac{n+2}{n+1}$$

as this is a Frullani integral. After changing the index $p \to p+1$ we get the following result.

**Corollary 4.** *For every integer $k \geq 0$ we have the series identities*

(18) $$\sum_{p=1}^{\infty} \frac{(-1)^{p-1} c_{p-1}}{p! p^k} = \sum_{p=0}^{\infty} \frac{(-1)^p c_p}{p!(p+1)^{k+1}} = \sum_{n=k}^{\infty} \begin{bmatrix} n \\ k \end{bmatrix} \frac{1}{n!} \ln\frac{n+2}{n+1}.$$

(19) $$\sum_{p=1}^{\infty} \frac{d_{p-1}}{p! p^k} = \sum_{p=0}^{\infty} \frac{d_p}{p!(p+1)^{k+1}} = \sum_{n=k}^{\infty} \begin{bmatrix} n \\ k \end{bmatrix} \frac{1}{n!} \ln\left(1+\frac{1}{n}\right).$$

*Proof.* The first identity has been proved above. For the second one we use the function

$$f(x) = \frac{-x}{(1-x)\ln(1-x)}$$

in the same way, qed.

Both series on the right hand sides of (18) and (19) are very slowly convergent series with positive terms. In (19), for instance, with $k \geq 1$ fixed

$$\begin{bmatrix} n \\ k \end{bmatrix} \frac{1}{n!} \ln\left(1+\frac{1}{n}\right) \sim \frac{1}{(k-1)!} \frac{(\gamma + \ln n)^{k-1}}{n} \ln\left(1+\frac{1}{n}\right).$$

### 3.4. Dirichlet series with derangement numbers

The derangement numbers

$$D_n = n! \sum_{j=0}^{n} \frac{(-1)^j}{j!}$$

are popular in combinatorics [7, p.180], [10, pp. 194-200]. We will relate them to Stirling numbers of the first kind. The generating function for the derangement numbers is given by

$$D(x) = \frac{e^{-x}}{1-x} = \sum_{n=0}^{\infty} D_n \frac{x^n}{n!} \quad (|x|<1).$$

In this case

$$\int_0^1 D(x)(1-x)^n \, dx = \int_0^1 e^{-x}(1-x)^{n-1} dx = \frac{1}{e}\int_0^1 e^t t^{n-1} dt = (-1)^n (n-1)! \left( e^{-1} - \sum_{j=0}^{n-1} \frac{(-1)^j}{j!} \right)$$

and therefore, we come to the series identity below.

**Corollary 5**. *For every integer* $k \geq 0$

(20) $$\sum_{p=1}^{\infty} \frac{D_{p-1}}{(p-1)! \, p^{k+1}} = \sum_{p=0}^{\infty} \frac{D_p}{p!(p+1)^{k+1}} = \sum_{n=k}^{\infty} \begin{bmatrix} n \\ k \end{bmatrix} \frac{(-1)^n}{n} \left( e^{-1} - \sum_{j=0}^{n-1} \frac{(-1)^j}{j!} \right).$$

From Taylor's formula we get the estimate

$$\left| e^{-1} - \sum_{j=0}^{n-1} \frac{(-1)^j}{j!} \right| \leq \frac{1}{n!}$$

which assures convergence for the last series in view of (10).

The computation of the above integral shows that with $f(x) = e^{-x}$ in (7) we come to

$$\int_0^1 e^{-x}(1-x)^n \, dx = (-1)^{n+1} n! \left( e^{-1} - \sum_{j=0}^{n} \frac{(-1)^j}{j!} \right)$$

and this result implies the identity

(21) $$\sum_{p=0}^{\infty} \frac{(-1)^p}{p!(p+1)^{k+1}} = \sum_{n=k}^{\infty} \begin{bmatrix} n \\ k \end{bmatrix} (-1)^{n+1} \left( e^{-1} - \sum_{j=0}^{n} \frac{(-1)^j}{j!} \right).$$

In the same way

$$\int_0^1 e^x (1-x)^n \, dx = n! \left( e - \sum_{j=0}^{n} \frac{1}{j!} \right)$$

and therefore,

(22) $$\sum_{p=0}^{\infty} \frac{1}{p!(p+1)^{k+1}} = \sum_{n=k}^{\infty} \begin{bmatrix} n \\ k \end{bmatrix} \left( e - \sum_{j=0}^{n} \frac{1}{j!} \right).$$

### 3.5. Identities for Dirichlet series with binomial and central binomial coefficients

Noticing that

$$\int_0^1 (1-x)^\lambda dx = \frac{1}{\lambda+1} \quad (\lambda > -1)$$

we take the function

$$f(x) = (1-x)^r = \sum_{n=0}^\infty \binom{r}{n}(-1)^n x^n, \quad |x|<1, \quad a_n = \binom{r}{n}(-1)^n$$

and from (7) we come to the next result:

**Corollary 6.** *For every integer $k \geq 0$ and $k+r+1>0$*

(23) $$\sum_{p=0}^\infty \binom{r}{p}\frac{(-1)^p}{(p+1)^{k+1}} = \sum_{n=k}^\infty \begin{bmatrix} n \\ k \end{bmatrix}\frac{1}{n!(n+r+1)}.$$

Here we need $n+r+1>0$. This will be true when $k+r+1>0$ or $r>-k-1$. This identity was obtained in [2, Example 10] by other means.

When $r$ is a nonnegative integer, the sum on the left hand side is finite. For example, when $r=0$ equation (23) turns into (15). For $r=1$ and $r=2$ we have correspondingly as in [2, Example 10]

(24) $$1 - \frac{1}{2^{k+1}} = \sum_{n=k}^\infty \begin{bmatrix} n \\ k \end{bmatrix}\frac{1}{n!(n+2)}$$

(25) $$1 - \frac{1}{2^k} + \frac{1}{3^{k+1}} = \sum_{n=k}^\infty \begin{bmatrix} n \\ k \end{bmatrix}\frac{1}{n!(n+3)}$$

etc.

For $r=-1/2$ and $r=1/2$ the binomial coefficients take a special form

$$\binom{-1/2}{p} = (-1)^p \binom{2p}{p}\frac{1}{4^p}, \quad \binom{1/2}{p} = (-1)^{p+1}\binom{2p}{p}\frac{1}{4^p(2p-1)}$$

and (23) produces the two identities involving central binomial coefficients:

**Corollary 7.** *For any $k \geq 0$*

(26) $$\sum_{p=0}^\infty \binom{2p}{p}\frac{1}{4^p(p+1)^{k+1}} = 2\sum_{n=k}^\infty \begin{bmatrix} n \\ k \end{bmatrix}\frac{1}{n!(2n+1)}$$

(27) $$\sum_{p=0}^\infty \binom{2p}{p}\frac{1}{4^p(2p-1)(p+1)^{k+1}} = -2\sum_{n=k}^\infty \begin{bmatrix} n \\ k \end{bmatrix}\frac{1}{n!(2n+3)}.$$

In particular, with $k=0$ in (26) we have the evaluation

$$\sum_{p=0}^{\infty} \binom{2p}{p} \frac{1}{4^p (p+1)} = 2$$

coming from the generating function of the Catalan numbers for $x = \frac{1}{4}$ (see next subsection).

Wang and Xu [17] studied series similar to the one on the left hand side in (25) and evaluated them in terms of multiple zeta values.

### 3.6. Dirichlet series with Catalan numbers

We involve now the Catalan numbers

$$C_p = \binom{2p}{p} \frac{1}{1+p}$$

which are very popular in combinatorics and analysis [7, p.101], [8, p.53], [10, 203], [4]. It is easy to see that for $p \geq 1$

$$\binom{2p}{p} \frac{1}{2p-1} = 2 C_{p-1}.$$

Our first identity with Catalan numbers comes from (27) written in the form

(28) $$\sum_{p=1}^{\infty} \frac{C_{p-1}}{4^p (p+1)^{k+1}} = \frac{1}{2} - \sum_{n=k}^{\infty} \begin{bmatrix} n \\ k \end{bmatrix} \frac{1}{n!(2n+3)}.$$

Correspondingly, for $k = 0, 1, 2$ we have from (28)

(29) $$\sum_{p=1}^{\infty} \frac{C_{p-1}}{4^p (p+1)} = \frac{1}{6}$$

(30) $$\sum_{p=1}^{\infty} \frac{C_{p-1}}{4^p (p+1)^2} = \frac{1}{2} - \sum_{n=1}^{\infty} \frac{1}{n(2n+3)} = \frac{2 \ln 2}{3} - \frac{7}{18}$$

(31) $$\sum_{p=1}^{\infty} \frac{C_{p-1}}{4^p (p+1)^3} = \frac{1}{2} - \sum_{n=2}^{\infty} \frac{H_{n-1}}{n(2n+3)} = -\frac{77}{54} + \frac{\pi^2}{18} + \frac{16}{9} \ln 2 - \frac{2}{3} \ln^2 2$$

where in equation (30) we used the evaluation (for example, from Wolfram Alpha)

$$\sum_{n=1}^{\infty} \frac{1}{n(2n+3)} = \frac{8}{9} - \frac{2 \ln 2}{3}.$$

and the second equality in (31) was found by Mathematica and was provided by one of the referees.

The generating function for the Catalan numbers is [4. 13]

$$\frac{2}{1+\sqrt{1-4x}} = \sum_{n=0}^{\infty} C_n x^n \quad (|x|<4)$$

Replacing $x$ by $x/4$ we consider the function

$$f(x) = \frac{2}{1+\sqrt{1-x}} = \sum_{n=0}^{\infty} \frac{C_n}{4^n} x^n$$

and apply our theorem to it. Using the substitution $1-x = t^2$ we find

$$\int_0^1 f(x)(1-x)^n dx = 2\int_0^1 \frac{(1-x)^n}{1+\sqrt{1-x}} dx = 4\int_0^1 \frac{t^{2n+1}}{1+t} dt = 4\beta(2n+2)$$

where

$$\beta(x) = \int_0^1 \frac{t^{x-1}}{1+t} dt = \sum_{m=0}^{\infty} \frac{(-1)^m}{m+x}$$

is Nielsen's beta function. We come to the curious companion to (28)

(32) $$\sum_{p=0}^{\infty} \frac{C_p}{4^p (p+1)^{k+1}} = 4 \sum_{n=k}^{\infty} \begin{bmatrix} n \\ k \end{bmatrix} \frac{\beta(2n+2)}{n!}.$$

For $k=0$ this gives the known identity

(33) $$\sum_{p=0}^{\infty} \frac{C_p}{4^p (p+1)} = 4\beta(2) = 4(1-\ln 2).$$

For $k=1$ in (32) we find

(34) $$\sum_{p=0}^{\infty} \frac{C_p}{4^p (p+1)^2} = 4 \sum_{n=1}^{\infty} \frac{\beta(2n+2)}{n}$$

and for $k=2$ we have a series identity involving Catalan, harmonic numbers, and beta values

(35) $$\sum_{p=0}^{\infty} \frac{C_p}{4^p (p+1)^3} = 4 \sum_{n=1}^{\infty} \frac{H_{n-1} \beta(2n+2)}{n}.$$

At the same time from (8)

(36) $$\sum_{p=0}^{\infty} \frac{C_p}{4^p (p+1)^{k+1}} = \frac{2(-1)^k}{k!} \int_0^1 \frac{(\ln x)^k}{1+\sqrt{1-x}} dx$$

and for $k=0$ this confirms (33). For $k=1$ by computing the integral we find

(37) $$\sum_{p=0}^{\infty}\frac{C_p}{4^p(p+1)^2}=8-8\ln 2+4(\ln 2)^2-\frac{\pi^2}{3}$$

which gives also the value of the series on the right hand side in (34)

$$\sum_{n=1}^{\infty}\frac{\beta(2n+2)}{n}=2-2\ln 2+(\ln 2)^2-\frac{\pi^2}{12}.$$

The integral $\int_0^1 \frac{\ln x}{1+\sqrt{1-x}}dx$ can be computed by the substitution $1-x=t^2$ followed by the expansion of $\ln(1-t^2)$ in power series. More directly, one can use Maple or Mathematica.

### 3.7. Dirichlet series with even central binomial coefficients

The numbers $\binom{4n}{2n}$ appear in some interesting applications in mathematics [5, 6]. Their generating function for $0\leq x<1$ is (see [6])

$$f(x)=\sum_{n=0}^{\infty}\binom{4n}{2n}\frac{x^n}{16^n}=\frac{1}{2}\left(\frac{1}{\sqrt{1-\sqrt{x}}}+\frac{1}{\sqrt{1+\sqrt{x}}}\right)=\frac{1}{2}\left(\frac{\sqrt{1-\sqrt{x}}+\sqrt{1+\sqrt{x}}}{\sqrt{1-x}}\right)$$

which is easily derived from the binomial series. The theorem implies the identity below.

**Corollary 8.** *For every integer $k\geq 0$ we have the identity*

(38) $$\sum_{p=0}^{\infty}\binom{4p}{2p}\frac{1}{16^p(p+1)^{k+1}}=\sum_{n=k}^{\infty}\begin{bmatrix}n\\k\end{bmatrix}\frac{A_n}{n!}$$

*where*

$$A_n=\frac{1}{2}\int_0^1\left(\sqrt{1-\sqrt{x}}+\sqrt{1+\sqrt{x}}\right)(1-x)^{n-\frac{1}{2}}dx.$$

This interesting integral supposedly has the form $A_n=a_n-b_n\sqrt{2}$ with $a_n$, $b_n$ positive rational numbers. This hypothesis was suggested by several cases verified by Maple.

A rough estimate gives

$$A_n\leq\frac{1}{2}\int_0^1(1+2)(1-x)^{n-\frac{1}{2}}dx=\frac{3}{2n+1}$$

which in view of (10) provides good convergence for the right hand side in (38). Also, from (8) we find the integral representation

$$\text{(39)} \quad \sum_{p=0}^{\infty}\binom{4p}{2p}\frac{1}{16^p(p+1)^{k+1}} = \frac{(-1)^k}{2k!}\int_0^1\left(\frac{1}{\sqrt{1-\sqrt{x}}} + \frac{1}{\sqrt{1+\sqrt{x}}}\right)(\ln x)^k\,dx$$

$$= \frac{(-1)^k 2^k}{k!}\int_0^1\left(\frac{t}{\sqrt{1-t}} + \frac{t}{\sqrt{1+t}}\right)(\ln t)^k\,dt$$

for every $k \geq 0$. For example, when $k = 0, 1$

$$\text{(40)} \quad \sum_{p=0}^{\infty}\binom{4p}{2p}\frac{1}{16^p(p+1)} = \frac{8}{3} - \frac{2}{3}\sqrt{2}$$

$$\text{(41)} \quad \sum_{p=0}^{\infty}\binom{4p}{2p}\frac{1}{16^p(p+1)^2} = \frac{80}{9} - \frac{32}{9}(\ln 8 + \sqrt{2}) + \frac{16}{3}\ln(\sqrt{2}+1)$$

etc. The integral $\int_0^1\left(\frac{t}{\sqrt{1-t}} + \frac{t}{\sqrt{1+t}}\right)\ln t\,dt$ for (41) was computed here by using both Maple and Mathematica.

In conclusion the author wants to express his deep gratitude to the referees for a number of valuable comments and suggestions that helped to improve the paper.

# References

[1] **Victor Adamchik,** On Stirling numbers and Euler sums, *J. Comput. Applied Math.*, 79 (1997), 119-130.

[2] **Khristo N. Boyadzhiev**, Stirling numbers and inverse factorial series, arXiv:2012.14546v1[math.NT].

[3] **Khristo N. Boyadzhiev**, New series identities with Cauchy, Stirling, and harmonic numbers, and Laguerre polynomials, *J. Integer Seq.*, 23 (2020), 20.11.7.

[4] **Khristo N. Boyadzhiev**, Series with Central Binomial Coefficients, Catalan Numbers, and Harmonic Numbers, *J. Integer Seq.*, 15 (2) (2012), 12.1.7.

[5] **John M. Campbell, Jacopo D'Aurizio and Jonathan Sondow**, Hypergeometry of the Parbelos, *Amer. Math. Montly*, 127 (2020), 23-32.

[6] **John M. Campbell, Jacopo D'Aurizio and Jonathan Sondow,** On the interplay among hypergeometric functions, complete elliptic integrals and Fourier-Legendre series expansions, *J. Math. Anal. Appl.*, 479(1) (2019), 90-121.

[7] **Louis Comtet,** *Advanced Combinatorics*, Kluwer, 1974.

[8] **John H. Conway, Richard Guy**, *The Book of Numbers,* Corrected edition, Copernicus, 1995.

[9] **Henry W. Gould, Temba Shonhiwa**, A catalog of interesting Dirichlet series, *Missouri J. Math. Sci.* 20 (1) (2008).

[10] **Ronald L. Graham, Donald E. Knuth, Oren Patashnik**, *Concrete Mathematics*, Addison-Wesley Publ. Co., New York, 1994.

[11] **Godfrey H. Hardy, Marcel Riesz**, *The general theory of Dirichlet's series.* Cambridge University Press, 1915.

[12] **Charles Jordan**, *Calculus of finite differences*, Chelsea, New York, 1950 (First edition Budapest 1939).

[13] **Derrick Henry Lehmer**, Interesting series involving the central binomial coefficient, *Amer. Math. Monthly,* 92 (1985), 449–457.

[14] **Donatella Merlini, Renzo Sprugnoli, M. Cecilia Verri**, The Cauchy numbers, *Discrete Math.*, 306 (2006), 1906-1920.

[15] **Ian Tweddle,** *James Stirling's Methodus Differentialis: An Annotated Translation of Stirling's Text*, Springer, New York, 2003.

[16] **Weiping Wang, Yao Chen**, Explicit formulas of sums involving harmonic numbers and Stirling numbers, *J. Difference Equ. Appl.*,26(2020), 1369-1397.

[17] **Weiping Wang, Ce Xu**, Alternating multiple zeta values, and explicit formulas of some Euler–Apéry-type series, *Eur. J. Comb.*,93, (March 2021), 103283.